\documentclass{article}

\usepackage{amsmath}
\usepackage{amssymb}
\usepackage{epsfig}

\usepackage{psfrag}
\newtheorem{theo}{Theorem}[section]

\newtheorem{lemma}[theo]{Lemma}
\newtheorem{theorem}[theo]{Theorem}

\newtheorem{rem}{Remark}
\newtheorem{algo}{Algorithm}
\newcommand{\N}{\mathbb{N}}

\newcommand{\R}{\mathbb{R}}
\newcommand{\C}{\mathbb{C}}

\newcommand{\para}{\parallel}

\newcommand{\dt}{ \frac{\partial}{\partial t}}

\newcommand{\beq}{\begin{equation}}
\newcommand{\eeq}{\end{equation}}
\newcommand{\beqn}{\begin{eqnarray}}
\newcommand{\eeqn}{\end{eqnarray}}

\font\QEDlogofont=msam10 at 10pt
\def\QEDlogo{\hbox{\QEDlogofont\char'003}}
\def\QEDblogo{\hbox{\QEDlogofont\char'004}}
\newif\ifnologo\nologofalse
\newif\iflogo
\newif\ifblogo\blogofalse
\newif\iftopprhead\topprheadfalse
\def\prooffont{\normalsize}
\newenvironment{proof}{\par\addvspace{6pt plus2pt}
\par
\noindent\prooffont{\it Proof}\hskip6pt\ignorespaces}{%
   \ifblogo\hskip1.2pt
            \QEDblogo
   \else
   \ifnologo
   \else
   \hfill
            \QEDlogo
   \fi\fi
\par\addvspace{6pt plus2pt}\global\topprheadfalse}%

\newcommand{\Hspace}{\mathcal{H}}
\newcommand{\Vspace}{\mathbf{V}}

\author{Julien Salomon,  Gabriel Turinici
\\
\\\vspace{6pt}{\em{ CEREMADE, UMR CNRS 7534, Universit\'e Paris IX,}}
\\\vspace{6pt}{\em{ Place du Marechal de Lattre de Tassigny,}}
\\\vspace{6pt}{\em{ 75775 Paris Cedex 16, France}} }

\begin{document}
\title{A monotonic method for solving nonlinear optimal control problems with concave dependence on the state}

\maketitle

\begin{abstract}
Initially introduced in the framework of quantum control, the
so-called {\it monotonic algorithms} have demonstrated excellent numerical performance
when dealing with bilinear optimal control problems. 
This paper presents a unified formulation that can be applied to
more nonlinear settings compatible with the hypothesis detailed below.
In this framework, we show that the well-posedness of
the general algorithm is related to  a nonlinear evolution equation.
We prove the existence of the solution to this equation 
and give important properties of the optimal control functional. Finally we show how the algorithm works for selected 
models from the literature and compare it with the gradient algorithm.
\end{abstract}



\section{Introduction}

This paper presents a general unified formulation of several algorithms
that were proposed in different areas of nonlinear control (see works in~\cite{kara,zr98, aimemfg,cdcCarlierSalomon}). 
Given a cost functional to minimize $J(v)$ depending on the control $v$, 
and a system described by a state function $X(t)$ solution of the evolution equation~\eqref{direct1} below,
these algorithms are iterative 
procedures that construct a sequence of solution candidates $v^{k}$ with the
important monotonic behavior $J(v^{k+1}) \le J(v^{k})$ ; the algorithms have 
been named after this property as "monotonic".
A convenient advantage of these procedures is that the monotonicity does not requires any 
additional computational effort, but results from the definition of the procedure itself. 

These algorithms have first been used in the field of quantum control
 where the dynamics is controlled by a laser field.
In this framework the function that associates to a control $v$
associates the final state $X(T)$ of the system is highly nonlinear.  
This induces poor performance of standard, gradient-based algorithms. The "monotonic
schemes" introduced in~\cite{hbref39,kara,zr98} were found to perform
excellently in this  setting.
These schemes were used in bi-linear situations i.e., when the operator $A(t,v(t))$ is linear in $v(t)$ and for a cost functional
$J(v) = G(X(T)) + C(v)$ that is sum of a part $G(X(T))$ quadratic in the final state $X(T)$ and a part $C(v)$ quadratic in the control $v$.
These were soon followed by 
variants as those in~\cite{hbref74,hbref75,schirmer,sugawara-JCP03,fujimura-sugawara-93} that
included situations where $G(X)$ has negative semi-definite Hessian but $C(v)$ was still quadratic in the
control and, most importantly, $A(t,v)$ was linear in $v$.

In works cited up to now, the function $G(X(T))$ depends only on the final state $X(T)$ 
but adaptations were proposed in~\cite{yukioshi,ohtsuki2} to deal
with the case where $G$ depends on the whole dynamics of the control process $X(t)$ at intermediary times 
or when the dynamics involve bilinear integro-differential equations.

Similar procedures were also proposed in different control
applications where the evolution equation is of parabolic type (see~\cite{cdcCarlierSalomon,aimemfg}) 
or mixed hyperbolic-parabolic (see~\cite{hbref76,ohtsuki:661}). 

Up to this point all works presented above considered bilinear situations i.e., 
the evolution equation is linear in the state $X(t)$ and $A(t,v)$
is linear in the control $v$; only recently different cases were documented in the literature 
where $A(t,v)$ is polynomial in the control $v$ up to power $3$; 
in~\cite{claudedionjcp,sugny1} specific monotonic procedures were proposed that were showed to work in this setting. 

A situation when $A(t,v)$ depends polynomially on a one-dimensional control $v(t) \in \R$ was proposed in~\cite{ohtsukinonlin}.
A model where the system is a nonlinear Bose-Einstein condensate was given
in~\cite{tannor2002}.

This paper continues in this direction and 
treats situations with arbitrary nonlinear $A(t,v)$ but still keeps the requirement that the explicit dependence of $J$ on $X$ be concave 
(see the hypothesis detailed below).
In all situations where monotonic algorithms were introduced the well-posedness 
of the algorithms
(i.e., the existence of $v^{k+1}$ given $v^k$)
was proved by ad-hoc techniques 
although the algebraic computations share similar points.
The purpose of this paper is to identify and exploit the similarities 
present in all these situations, and
present a general setting that includes
the "monotonic" algorithms. This allows to tackle a
large class of non-linear situations that cannot be solved with techniques present 
today in the literature. We prove rigorously the existence and convergence of a 
procedure that from a control $v^k$ constructs a control $v^{k+1}$ such that the cost functional is monotonic. 
The question of whether such a procedure exists has never been asked before in the literature 
because up to now the authors considered only particular cost functionals $J$ and particular evolution equations; 
in each case they proposed explicit analytic formulaes for $v^{k+1}$ applicable to their situation. On the contrary, 
we show here that in all contexts covered by the theory (i.e., satisfy the hypothesis below) a control $v^{k+1}$ can always be found to ensure
$J(v^{k+1}) \le J(v^{k})$ and we also give a constructive procedure to
compute it.

The paper is structured as follows: Section~\ref{Sec:setting} defines
the general framework in which where our procedure applies. The algorithm is
presented in Section~\ref{Sec:MS}. At this point we show that the well-posedness of
the algorithm is related to  a nonlinear evolution equation
and prove the existence of the solution to this equation. We also give important properties of the
optimal control functional. 
Some examples of concrete realizations follow in Section~\ref{Sec:examples} together with numerical results illustrating
the application and the efficiency of the algorithm.

\section{Problem formulation}\label{Sec:setting}

Let $E$, $\Hspace$ and $\Vspace$ be Hilbert spaces with $\Vspace$ densely included in $\Hspace$.
We denote by $\cdot_E $ and $\langle \cdot, \cdot \rangle_{\Vspace}$
the scalar product associated with $E$ and $\Vspace$.

For any vector spaces ${\mathcal
  A}$ and ${\mathcal B}$, we denote by ${\mathcal L}({\mathcal A},{\mathcal B})$ the space
of linear continuous operators between  ${\mathcal
  A}$ and ${\mathcal B}$.

Given a real valued function $\varphi$, we denote by $\nabla_x\varphi$ its
  gradient with respect to the variable $x$. We also denote by
  $D_x$ and $D_{x,x}$ the first and the second derivative of
  vectorial functionals in the Fr\'echet sense.
\begin{rem}
Recall that, given $H_1$ and $H_2$ two Banach spaces and $U\subset
H_1$ an open subset of $H_1$, a function $f : U \to H_2$  is said to be 
Fr\'echet differentiable at $x \in U$ if there exists a continuous linear
operator $A_x \in {\mathcal L}(H_1, H_2)$ such that 
$$    \lim_{h \to 0} \frac{ \| f(x + h) - f(x) - A_x (h) \|_{H_2} }{ \|h\|_{H_1} } = 0. $$
The operator $A_x$ is then called the Fr\'echet differential (or Fr\'echet derivative) of $f$ at $x$
and is denoted $D_x f \triangleq A_x$.

Let us also recall that given an open set $\Omega \subset \R^\gamma$ and a Hilbert space $H_1$,  
the set $L^\infty(\Omega;H_1)$ is the space of functions $f$ from $\Omega$ with values 
in the Hilbert space $H_1$ such that for almost all $t \in \Omega$ the norm  $\| f(t) \|_{H_1}$ is bounded by
the same constant (the lowest of which is the $L^\infty(\Omega;H_1)$ norm of $f$).
One can likewise define $L^2(\Omega;H_1)$: 
\beq 
L^2(\Omega;H_1)= \{ f:\Omega\to H_1 \textrm{ such that } \int_\Omega \|f(t) \|^2_{H_1} dt < \infty \}.
\eeq 

When the derivatives of $f$ are considered the Sobolev spaces $W^{1,\infty}$ have to be introduced; we refer to~\cite{adams,yosida} for further details.
\end{rem}

Within an optimal control formulation, the control of a 
system described by a state function $X(t)$ is encoded in the following
optimization problem:
\begin{equation}\label{initpb}
\min_v J(v), 
\end{equation} 
where
\beq\label{eqn:F}
J(v) \triangleq \int_0^T F\big(t,v(t),X(t)\big)dt+G\big(X(T)\big).
\eeq
The functions $F:\R \times E \times \Vspace \rightarrow \R$ and 
$G:\Vspace \rightarrow \R$ are assumed to be differentiable and
integral assumed to exist. The state function $X(t) \in
\Vspace$ satisfies the following evolution equation
\begin{eqnarray}
\partial_t{X}+A(t,v(t))X&=&B(t,v(t))\label{direct1}\\
X(0)=X_0.\label{direct2}
\end{eqnarray}
where $v:[0,T]\to E$ is the control.
The possibly unbounded operator
$A(t,v): \R \times E \times \Hspace \to \Hspace$
is such that for almost all $t\in [0,T]$ the domain 
of $A(t,v)^{1/2}$ includes $\Vspace$; furthermore we take 
$B(t,v)$ such that for almost all $t\in [0,T]$ and all $v \in E$ 
we have  $B(t,v) \in {\mathcal L}(\Hspace,\Hspace) \cap {\mathcal
  L}(\Vspace,\Vspace^*)$\footnote{For a space $\Vspace$ we denote by
  $\Vspace^*$ its dual space.}. We postpone to Section~\ref{Sec:MS}
(cf. Lemma~\ref{lemma:existence4}, Theorem~\ref{thm:existence}) the
precise formulation of   
additional regularity assumptions to be imposed on $A,B,F,G$.

\begin{rem}
Note that $E$ is not necessarily of dimension one, not even finite dimensional, cf. Section~\ref{sec:exmfg}. This means
that the control can be a set of {\bf several} time-dependent functions but also a distributed control
depending on a spacial variable.
\end{rem}

 Let us stress that although the equation is linear in $X$ (for $v$ fixed) 
the mapping $v\mapsto X$ is {\bf not} linear ; the term $A(t,v(t))$ {\bf multiplies} 
the state $X$ and as such the mapping is highly nonlinear (of non-commuting exponential type).

\begin{rem} \label{rem:nonbilinear}
Most of the previous works considered a bilinear operator $A(t,v)$ i.e., $A(t,v)X=vX$; the only exceptions
(cf. discussion in the Introduction) were of the polynomial type (of order at most $3$ in~\cite{claudedionjcp,sugny1} and polynomial 
with $E=\R^1$ in~\cite{ohtsukinonlin}). The techniques present in the above papers
cannot be used for general operators $A(t,v)$ and control sets $E$. On
the contrary the results in this work  
include all the situations considered in the bibliography but also 
apply to nonlinearities in $v$ compatible with the hypothesis of
Lemma~\ref{lemma:existence4} and Thm.~\ref{thm:existence} below. 
\end{rem}

The following concavity with respect to $X$ will be assumed throughout the paper:
\beqn
& \ & \label{concavG}
\forall X,X'\in \Vspace,\ G(X') - G(X) \leq \langle \nabla_X G (X),
X'-X\rangle_{\Vspace} ,
\\ & \ & \label{concavF}
\forall t\in \R, \forall v\in E, \forall X,X'\in \Vspace,\nonumber
\\ &\ & F(t,v,X')-F(t,v,X) \le \langle \nabla_X F
(t,v,X), X'-X\rangle_{\Vspace}.
\eeqn
\begin{rem} \label{rem:hyp}
Unlike the more technical hypothesis that will be assumed latter, the properties 
\eqref{concavG},
\eqref{concavF} and the linearity of \eqref{direct1} are crucial to the existence of the {\it monotonic} algorithms. 
\end{rem}
\section{Monotonic algorithms}\label{Sec:MS}
We now present the structure of our optimization procedure together with
the general algorithm.

\subsection{Tools for monotonic algorithms}
The monotonic algorithms exploit a specific factorization which is the consequence
of the results described in this
section. To ease the notations, we will make explicit the dependence of $X$ on $v$, i.e. we will
write $X_v$ instead of $X$ in Eqs.~(\ref{direct1}--\ref{direct2}).

We define the adjoint state $Y_v$~(see \cite{gabasovbook,lionsedp})
which is the solution of the following evolution equation:
\begin{eqnarray} 
\partial_t{Y_v}-A^*\big(t,v(t)\big)Y_v+\nabla_{X}F\big(t,v(t),X_v(t)\big)&=&0\label{retrograde1}\\
Y_v(T)=\nabla_X G\big(X_v(T)\big).\label{retrograde2}
\end{eqnarray}
A first estimate about the variations in $J$
can be obtained: 

\begin{lemma}\label{mainlem}
For any $v',v:[0,T]\rightarrow E$ denote
\beqn \label{eq:delta00}
& \ & 
\Upsilon \Big(t,X_v(t),v(t),v'(t),Y_v(t),X_{v'}(t) \Big) \triangleq \nonumber\\
& \ & -\langle Y_v(t) ,  \Big(A\big(t,v'(t)\big)-A\big(t,v(t)\big)\Big)X_{v'}(t)\rangle_{\Vspace} \nonumber\\ & \ &
 \langle Y_v(t),B\big(t,v'(t)\big)-B\big(t,v(t)\big)\rangle_{\Vspace}\\
& \ & +F\big(t,v'(t),X_{v'}(t)\big)-F\big(t,v(t),X_{v'}(t)\big).
\eeqn
Then
\beq \label{eq:delta0}
   J(v')-J(v) \le     \int_0^T 
\Upsilon \Big( t,X_v(t),v(t),v'(t),Y_v(t),X_{v'}(t) \Big) dt.
\eeq
\end{lemma}
\proof Using successively  \eqref{concavG},\eqref{concavF}, \eqref{direct1} and finally \eqref{retrograde2}, we find that:
\beqn
J(v')-J(v)           & =  &  \int_0^T F\big(t,v(t),X_{v'}(t)\big)-F\big(t,v(t),X_{v}(t)\big)\nonumber\\&&\phantom{\int_0^T} + F\big(t,v'(t),X_{v'}(t)\big)-F\big(t,v(t),X_{v'}(t)\big)dt\nonumber\\
\phantom{J(v')-J(v)} &    & + G\big(X_{v'}(T)\big)-G\big(X_v(T)\big) \nonumber\\
\phantom{J(v')-J(v)} &\leq&  \int_0^T \!\!\!\langle \nabla_X F\big(t,v(t),X_v(t)\big), X_{v'}(t)-X_v(t)\rangle_{\Vspace}\nonumber\\
                     &    &  \phantom{ \int_0^T}+ F\big(t,v'(t),X_{v'}(t)\big)-F\big(t,v(t),X_{v'}(t)\big) dt \nonumber\\
\phantom{J(v')-J(v)} &    & +\langle Y_v(T),X_{v'}(T)-X_v(T)\rangle_{\Vspace}\nonumber\\
\phantom{J(v')-J(v)} &\leq&  \int_0^T\!\!\! \langle
\dt{Y_v}(t)-A\big(t,v(t)\big)^*Y_v(t) \nonumber \\
\phantom{J(v')-J(v)} &    & \phantom{\int_0^T}\!\!\! + \nabla_X F\big(t,v(t),X_v(t)\big),X_{v'}(t)-X_v(t)\rangle_{\Vspace} \nonumber\\
\phantom{J(v')-J(v)} &    & \phantom{\int_0^T}\!\!\! - \langle Y_v(t),\Big(A\big(t,v'(t)\big)-A\big(t,v(t)\big)\Big)X_{v'}(t)\rangle_{\Vspace} \nonumber \\
\phantom{J(v')-J(v)} &    & \phantom{\int_0^T}\!\!\! + \langle
Y_v(t),B(t,v'(t))-B(t,v(t)) \rangle_{\Vspace} \nonumber \\
                     &    & \phantom{\int_0^T}\!\!\! + F\big(t,v'(t),X_{v'}(t)\big)-F\big(t,v(t),X_{v'}(t)\big)dt.\nonumber
\eeqn
Due to \eqref{retrograde1}, the first term of the right-hand side of
this last inequality cancels and the result follows.\endproof

\begin{rem} \label{rem:dirac}
The purpose of this result is not to obtain an estimation of the increment $J(v')-J(v)$ 
{\it via} the adjoint
(which is well documented in optimal control theory, cf.~\cite{gabasovbook,lionsedp}); 
we rather emphasis that the evaluation of the integrand $\Upsilon$ 
at time $t$ requires information on the control $v(s)$ for all $s \in [0,T]$ (in order to compute 
$X_v(T)$ then
$Y_v(t)$) but on the second control $v'(s)$ only for $s \in [0,t]$ (because this is enough to
compute $X_{v'}(t)$). 
This estimate can be useful to decide {\bf at time $t$}
if the current value of the control $v'(t)$ will imply an increase or decrease of $J(v')$.
This localization property is a consequence of the concavity of $F$
and $G$ (in $X$) and bi-linearity induced by $A$. The 
purpose of the paper is to construct and theoretically support a 
general numerical algorithm that exploits this remark.
\end{rem}

\begin{rem} We can intuitively note that 
$\Upsilon$ has the 
factorized form:
\beqn\label{eq:factorization}
\Upsilon \Big( t,X_v(t),v(t),v'(t),Y_v(t),X_{v'}(t) \Big) = \Delta(v,v')(t) \cdot_E \big(v'(t)-v(t)\big), 
\eeqn
with $\cdot_E$ the $E$ scalar product.
Thus $v'$ can always be chosen so as to make it negative 
(in the worse case set it null by the choice $v'=v$). 
We will come back with a formal definition of $\Delta(v,v')(t)$  and a proof of the previous relation in Section~\ref{sec:existence}.
\end{rem}

A more general formulation can be obtained if we suppose that the
backward propagation of the adjoint state is performed with the
intermediate field $\widetilde{v}$ (cf. also \cite{jcp}), i.e. according to the equation :
\begin{eqnarray}
\dt{Y}_{\widetilde{v}}-A^*\big(t,\widetilde{v}(t)\big)Y_{\widetilde{v}}+\nabla_{X}F\big(t,v(t),X_{v}(t)\big)&=&0\nonumber\\
Y_{\widetilde{v}}(T)=\nabla_X G\big(X_{v}(T)\big).\nonumber
\end{eqnarray}
Note that because of its final condition, $Y_{\widetilde{v}}$ actually
also depends on $v$. Nevertheless, for the sake of simplicity, we keep
the previous notation. We then obtain the following lemma whose proof we let as exercise for the reader:

\begin{lemma}\label{mainlem2}
For any $v',\widetilde{v},v:[0,T]\rightarrow E$, 
\beqn
  J(v')-J(v)&\leq&\int_0^T-\langle Y_{\widetilde{v}}(t),\Big(A\big(t,v'(t)\big)-A\big(t,\widetilde{v}(t)\big)\Big)X_{v'}(t)\rangle_{\Vspace} \nonumber\\
            &    &\phantom{\int_0^T}+\langle Y_{\widetilde{v}}(t),B\big(t,v'(t)\big)-B\big(t,\widetilde{v}(t)\big)\rangle_{\Vspace} \nonumber\\
            &    &\phantom{\int_0^T}+F\big(t,v'(t),X_{v'}(t)\big)-F\big(t,\widetilde{v}(t),X_{v'}(t)\big)dt\nonumber\\
            &    &+\int_0^T-\langle Y_{\widetilde{v}}(t),\Big(A\big(t,\widetilde{v}(t)\big)-A\big(t,v(t)\big)\Big)X_{v}(t)\rangle_{\Vspace} \nonumber\\
            &    &\phantom{+\int_0^T}+\langle Y_{\widetilde{v}(t)}(t),B\big(t,\widetilde{v}(t)\big)-B\big(t,v(t)\big)\rangle_{\Vspace}  \nonumber\\
            &    &\phantom{+\int_0^T}+F\big(t,\widetilde{v}(t),X_{\widetilde{v}}(t)\big)-F\big(t,v(t),X_{\widetilde{v}}(t)\big)dt.\nonumber
\eeqn
\end{lemma}
In this lemma, the variation in the cost functional $J$ is expressed
as the sum of two terms, and can be considered as factorized with
respect to $v'-\widetilde{v}$ and $\widetilde{v}-v$.

\begin{rem} \label{rem:lemmanews}
Lemmas~\ref{mainlem} and~\ref{mainlem2} are generalizations of previous results that were proved in the bilinear case.
To the best of our knowledge, only specific corollaries requiring 
additional assumptions have appeared in the literature up to now.
\end{rem}

\subsection{The algorithms} \label{sec:algo}
The factorization \eqref{eq:factorization} that will be proved in
Lemma~\ref{factorlemma} enables to design various ways to ensure that
$J(v')\leq 
J(v)$, i.e. that guaranty the monotonicity resulting from the update $v'\leftarrow
v$. This allows to present a general structure for our
class of optimization algorithms. We focus on the one that results 
from Lemma \ref{mainlem}.

\begin{algo}(Monotonic algorithm)\label{algoms}\\
Given an initial control $v^0$, the sequence $(v^{k})_{k\in\N}$ is
computed iteratively by:
\begin{enumerate}
\item Compute the solution $X_{v^k}$ of (\ref{direct1}--\ref{direct2}) with $v=v^k$.
\item Compute the solution $Y_{v^k}$ of (\ref{retrograde1}--\ref{retrograde2}) with $v=v^k$ backward in time 
from 
$$Y_{v^k}(T) \triangleq \nabla_X G\big(X_{v^k}(T)\big).$$
\item Define (as explained latter) $v^{k+1}$ together with $X_{v^{k+1}}$ such that for all $t\le
  T$ the following monotonicity condition be satisfied :
\beq\label{monotoniccond}
\Delta (v^{k+1},v^{k})(t) \cdot_E \Big(v^{k+1}(t)-v^k(t)\Big)
\le 0.
\eeq
\end{enumerate}
\end{algo}
\noindent Lemma \ref{mainlem} then guarantees that $J(v^{k+1})\leq J(v^k)$. Several strategies can be used to ensure
\eqref{monotoniccond}; we will present one below. Its importance stems from the fact that no further optimization is necessary once 
this condition is fulfilled. 
In order to guarantee \eqref{monotoniccond}, many authors
(see~\cite{jcp,kara,zr98}) consider an update formula of the form:
\begin{equation}\label{usualupdate}
 v^{k+1}(t)-v^k(t)=  - \frac{1}{\theta}  \Delta(v^{k+1},v^{k})(t),
\end{equation}
where $\theta$ is a positive number, that can also depend on $k$ and
$t$.  In what follows, we focus on the existence of solution of \eqref{usualupdate}, and on practical methods to compute it.
If $v^{k+1}$ satisfies \eqref{usualupdate}, the variations in $J$ satisfy:
\beqn
  J(v^{k+1})-J(v^k) & \leq &  -\theta  \int_0^T (v^{k+1}(t)-v^k(t))^2 dt.\nonumber
\eeqn
Note that \eqref{usualupdate} reads as an update formula
combining on the one hand a gradient method:
$$v^{k+1}(t)-v^k(t)=- { \frac{1}{\theta}} \Delta(v^{k}  ,v^{k} )(t), $$
and on the other hand the so-called Proximal Algorithm (as described in~\cite{bolte_attouch}) which prescribes:
$$v^{k+1}(t)-v^k(t)=- { \frac{1}{\theta}} \Delta(v^{k+1},v^{k+1})(t). $$
\begin{rem} When $F=0$ and $A$ is independent of $v$ i.e., linear control with
final objective, \eqref{usualupdate} coincides with a gradient method.
\end{rem}

\subsection{Well-posedness of the algorithm} \label{sec:existence}
In this section, we focus on the procedure obtained when using Algorithm
\ref{algoms} with the update
formula \eqref{usualupdate}.
To the best of our knowledge no theoretical result exists in the literature to prove the existence
of a solution to the Eq.~(\ref{usualupdate}) for general choices of $A(t,v)$ and general space of controls $E$ because previous authors only dealt with
particular choices of functionals $F,G$, operators $A,B$ and managed to obtain in each case an analytic solution;
we provide here such a proof together with a convergent procedure to compute it.
Since this procedure involves the resolution of an implicit
equation, the proof is non-trivial and has been split in three parts: two preparatory Lemmas (\ref{factorlemma} and \ref{lemma:existence4})
and the final result in Theorem~\ref{thm:existence}. As a by-product, we obtain a proof of the monotonicity of the algorithm.

\begin{lemma}\label{factorlemma}
Suppose that for any $t\in[0,T]$:

- ${\mathcal A}:\R \times \Vspace \times \Vspace \times E \to\R$
defined by ${\mathcal A}(t,X,Y,v) = \langle Y, A(t,v) X
\rangle_{\Vspace}$ is of $C^1$ class with respect to $v$ for any $X,Y,v$;

- ${\mathcal B}:\R \times \Vspace \times E\to\R$ with ${\mathcal
  B}(t,Y,v) = \langle Y, B(t,v) \rangle_{\Vspace}$ is of $C^1$ class
 with respect to $v$ for any $Y,v$;

- $F$ is of $C^1$ class with respect to $v\in E$ for any  $X,Y,v$. 

Then there
exists $\Delta(\cdot,\cdot;t,X,Y) \in C^0(E^2,E)$ 
such that, for all $v,v'\in E$ 
\beqn \label{eq:delta}
\Delta(v',v;t,X,Y) \cdot_E \Big(v'-v\Big) &=&
-\left\langle Y,
  \Big(A(t,v')-A(t,v) \Big) X+B(t,v')-B(t,v)\right\rangle_{\Vspace}  
\nonumber \\ & \ & +F(t,v',X)-F(t,v,X).
\eeqn
Moreover,
 $\Delta(\cdot,\cdot;t,X,Y)$ can be defined through the explicit formula:
\beqn 
\Delta(v',v;t,X,Y)& =& \int_0^1 
-\nabla_w\Big( \left\langle Y,
   A(t,w)  X -
   B(t,w)\right\rangle_{\Vspace}\Big)\Big|_{w=v+\lambda(v'-v)}   \nonumber 
\\& \ &
\phantom{- \int_0^1 }+\nabla_v F(t,v+\lambda(v'-v),X) d\lambda.\label{eq:diffdelta}
\eeqn
\end{lemma}

\proof We denote by $\|\cdot\|$ the norm associated with $E$.
Since $\mathcal{A,B}$,$F$ are Fr\'echet differentiable with respect to $v$ the full expression in Eq.~\eqref{eq:delta}
is of the form ${\Xi}(v')-{\Xi}(v)$ with  ${\Xi}(v) = -\mathcal{A}(t,X,Y,v) + \mathcal{B}(t,Y,v)  - F(t,v,X) $ differentiable with respect to
$v$; 
Eq.~\eqref{eq:diffdelta} is an application of the identity
$${\Xi}(v')-{\Xi}(v) = \int_0^1 \nabla_v{\Xi}(v+\lambda(v'-v) )
d\lambda\cdot_E (v'-v) .$$
The continuity is obtained from that of $\nabla_v {\Xi}$.
\endproof

\begin{lemma} \label{lemma:existence4}
Suppose that 

- $\mathcal{A,B},F$ are of (Fr\'echet) $C^2$ class with respect to $v$ with $D_{vv}A$, $D_{vv}B$ uniformly bounded 
as soon as $X$,$Y$ are in a bounded set;

- $\nabla_v F$ is of $C^1$ class with respect to $X$;

-  $D_{vv}F(t,\cdot,X)$ is bounded by a positive, continuous,
increasing function $X\mapsto k(\|X\|)$. 

Then given $\varepsilon>0$, $(t,v,X,Y)\in \R \times E \times \Vspace \times \Vspace$ and a
bounded neighborhood $W$ of $(t,v,X,Y)$, there exists $\theta^\star>
0$ depending only on $\varepsilon$, $W$, $\|v\|$, $\|X\|$ and
$\|Y\|$ such that, for any
$\theta > \theta^\star$
\begin{enumerate}
\item \label{item:existence}
$\Delta(v',v;t,X,Y) = - \theta  (v'-v)$ has a unique solution 
$v'={\cal V}_\theta(t,v,X,Y) \in E$;
\item ${\cal V}_\theta(t,v,X,Y)=v$ implies 
\beq \label{eq:critv}
 - \nabla_v\Big(\left\langle Y, A(t, v)  X\right\rangle_{\Vspace}\Big)(v) +
 \nabla_v\Big(\left\langle Y, B(t,v) \right\rangle_{\Vspace}\Big)(v)   
+\nabla_vF(t,v,X) =0,
\eeq

\item \label{item:estimation1}
$\| {\cal V}_\theta(t,v,X,Y) -v \| \le
\frac{\| X\| \| Y\|+\| Y\|+k(\|X\|)}{\theta}
\{ M_0(t) + M_1\| v \|\} $ with $M_0(t)$ and $M_1$ independent of $v,X,Y$.
 If the dependence of $A,B,F$ with respect to 
$t$ is smooth then $M_0(t)$ is bounded on $[0,T]$;

\item ${\cal V}_\theta(t,v,X,Y)$ is continuous on $W$;

\item \label{item:estimation2}
Let $X$ belong to a bounded set; then 
$X \mapsto {\cal V}_\theta(t,v,X,Y)$ is Lipschitz with the Lipschitz
constant smaller than $\varepsilon$.

\end{enumerate}
\end{lemma}

\proof 
\begin{enumerate}
\item
Denote $h = v'-v$ and define ${\cal G}_{t,v,X,Y}(h) =\frac{-\Delta(v+h,v;t,X,Y)}{ \theta 
} $. When the dependence is clear we will write simply
${\cal G}(h) $ instead of ${\cal G}_{t,v,X,Y}(h) $.
We look thus for a solution to the following fixed point problem:
${\cal G}(h) = h$. For $\theta$ 
large enough, the mapping ${\cal G}$ 
is a (strict) contraction and we obtain the conclusion
by a Picard iteration. The uniqueness is a consequence of the contractivity of ${\cal G}$.

\item
If $v'=v$ then $h=0$ thus ${\cal G}(h)=0$ which gives \eqref{eq:critv} after using \eqref{eq:diffdelta}.

\item For $\theta$ large enough, the mapping ${\cal G}$ is not only a contraction but
has its Lipschitz constant less than, say, $1/2$. 
Because of the contractivity of ${\cal G}$, we have
$ \| h \| - \|{\cal G}(0)\| \le \| h - {\cal G}(0)\| = \| {\cal G}(h) - {\cal G}(0)\| \le \frac{1}{2} \| h \| $,
which amounts to $\| h \| \le 2 \|{\cal G}(0)\|$. Next, we note that 
\beqn
& \ & 
 \|{\cal G}(0)\|  \le \frac{ \| \Delta(v,v,t,X,Y) - \Delta(0,0,t,X,Y)\| + 
\|\Delta(0,0,t,X,Y)\| }{\theta  
}
\nonumber \\ & \ &  \le M_2\| v \|+M_3(t)\nonumber
\eeqn
and the estimate follows.
\item Formula
\eqref{eq:diffdelta} shows that $\Delta$ depends continuously on
$t,v,v',X,Y$.  
Consider converging sequences $t_n \to t$, $v_n \to v$, $X_n \to X$, $Y_n \to Y$ 
  and
define $h_n \triangleq {\cal V}_\theta
(t_n,v_n,X_n,Y_n)$ and $h \triangleq {\cal V}_\theta(t,v,X,Y)$. \\
Given $W$ and $\eta>0$, consider a large enough value of $\theta$ such
that:\\

- for any $(t',v',X',Y')\in W$, ${\cal G}_{t',v',X',Y'}$ is a
contraction with Lipschitz constant less than 1/2.\\

- for any $(t',v',X',Y')$,  $(t'',v'',X'',Y'')\in W$,
$$\| \Delta(v''+ h ,v'',t'',X'',Y'')-\Delta(v'+h ,v',t',X',Y') \| \le
\eta.$$ 

This last property implies  
$\| {\cal G}_{t_n,v_n,X_n,Y_n}(h) - {\cal G}_{t,v,X,Y}(h) \| \le \frac{\eta}{\theta}$ for 
$n$ large enough.
On the other hand 
\beqn 
\|h_n - h \| &=& \| {\cal G}_{t_n,v_n,X_n,Y_n}(h_n) - {\cal G}_{t,v,X,Y}(h) \|  \nonumber\\
&\le& \| {\cal G}_{t_n,v_n,X_n,Y_n}(h_n) - {\cal G}_{t_n,v_n,X_n,Y_n}(h)\|\nonumber\\
&&+\| {\cal G}_{t_n,v_n,X_n,Y_n}(h) - {\cal G}_{t,v,X,Y}(h) \| \nonumber\\
&\le& \frac{1}{2} \|h_n - h \| +\frac{\eta}{\theta}
.\nonumber
\eeqn 
We have thus obtained that for $n$ large enough  :  
$\frac{1}{2} \|h_n - h \| \le \frac{\eta}{\theta}$ and the continuity follows.
\item
 Subtracting the two equalities
$$
\Delta(V_1,v;t,X_1,Y) = - \theta (V_1-v), \ 
\Delta(V_2,v;t,X_2,Y) = - \theta (V_2-v) 
$$
and using that 
$\Delta(V,v;t,X,Y)$ is $C^1$ in $X$ and $v$
gives to first order
$$
\Delta_V(...)(V_1-V_2) + \Delta_X(...)(X_1-X_2)=-\theta (V_1-V_2).
$$
For $\theta$  large enough the operator $\Delta_V(...) + \theta \cdot Id$ is invertible and
the conclusion follows.
\end{enumerate}

\endproof

\begin{rem} Note that $\theta^\star$ is proportional to $( \| X\|_\Vspace \| Y\|_\Vspace+\| Y\|_\Vspace+k(\|X\|_\Vspace))  $.
\end{rem}

We are now able to construct a procedure such that the existence of
$v^{k+1}(t)$ satisfying \eqref{monotoniccond} is guaranteed.

\begin{theorem} \label{thm:existence}
Suppose that $A,B,F$ satisfy the hypothesis of the Lemma~\ref{lemma:existence4}. Also suppose that 
the operators $A,B$ are such that Eqs. (\ref{direct1}--\ref{direct2}) and (\ref{retrograde1}-\ref{retrograde2})
have solutions for any $v\in L^\infty(0,T;E)$ with $v\mapsto X$, $v\mapsto Y$ 
locally Lipschitz. Then:
\begin{enumerate}
\item For any $v\in L^\infty(0,T;E)$, there exists $\theta^\star>0$ such that for any $\theta >
  \theta^\star$, the (nonlinear) evolution system
\begin{eqnarray} 
& \ & \partial_t{X_{v'}}(t)+A(t,v')X_{v'}(t)=B(t,v') \label{eq:Xnl1}\\
& \ & v'(t) = {\cal V}_\theta(t,v(t),X_{v'}(t),Y_v(t)) \label{eq:Xnl2}\\
& \ & X_{v'}(0)=X_0\label{eq:Xnl3}
\end{eqnarray}
has a solution. Here 
$Y_v$ is the adjoint state defined by
(\ref{retrograde1}--\ref{retrograde2}) and corresponding to control $v$.

\item There exists a sequence $(\theta_k)_{k\in\N}$ such that
the algorithm~\ref{algoms} (cf. Section~\ref{sec:algo})  

a/ initialization  $v^0 \in  L^\infty(0,T;E)$, 

b/  $v^{k+1}(t)={\cal V}_{\theta_k}(t,v^k(t),X_{v^{k+1}}(t),Y_{v^k}(t))$ 

is monotonic and satisfies

\beq\nonumber
J(v^{k+1})-J(v^k)  \leq  - \theta_k \|v^{k+1}-v^k \|_{L^2([0,T])}^2.
\eeq

\item With the notations above, if for all $t\in [0,T]$ $v^{k+1}(t)=v^k(t)$ (i.e. algorithm stops) then $v^k$ is a critical point of $J$:
$\nabla_{v}J(v^k) = 0$.

\end{enumerate}
\end{theorem}
\proof 
Some of the proof is contained in the previous lemmas. The part that still has to be
proved is the existence of a solution to~\eqref{eq:Xnl1}-\eqref{eq:Xnl3}.

Given $v\in L^\infty(0,T;E)$, consider the following iterative procedure : 
\beq \nonumber
v_0=v, \ v_{l+1}(t) = {\cal V}_\theta(t,v(t),X_{v_l}(t),Y_v(t)).
\eeq
We take a spherical neighborhood $B_v(R)$ of $v$ of radius $R$ and
suppose that  $\forall k \le l,\ v_k\in B_v(R)$.
Since the
correspondence $v \mapsto X_v$ is continuous, it follows that  
the set of solutions $S_{v,R} \triangleq \{ X_w; w \in B_v(R)\} $ of  \eqref{direct1}
is bounded.
\ In particular for $w=v_l$
by the item~\ref{item:estimation1} of Lemma~\ref{lemma:existence4}
the quantity $\|{\cal V}_\theta(t,v(t),X_{v_l}(t),Y_v(t)) - v \|$
will be bounded by $R$ for $\theta$  large enough (depending on $R$, independent of $l$),
i.e. $v_{l+1} \in B_v(R)$.
Thus $v_l \in B_v(R)$ for all $l \ge 1$.

Since  $S_{v,R}$ is bounded, 
recall that by item~\ref{item:estimation2} of Lemma~\ref{lemma:existence4}
the mapping 
$X \mapsto {\cal V}_\theta(t,v(t),X,Y_v(t))$ has on $S_{v,R}$ a
Lipschitz constant as small as desired. Since the mapping
$w \mapsto X_w$ is  Lipschitz, for $\theta$ large enough,
$w \in B_v(R) \mapsto {\cal V}_\theta(t,v(t),X_w,Y_v(t))$ is a contraction. By a Picard argument
the sequence $v_l$ is converging. The limit will
 be the solution of~(\ref{eq:Xnl1}--\ref{eq:Xnl2}).\endproof

\section{Examples}\label{Sec:examples} 
We now present three examples that fit into the setting of
Theorem~\ref{thm:existence}. The space does not allow to treat different variants (cf. references
in Introduction) so we leave them as an exercise to the reader. 

Within the framework of control theory, nonlinear formulations prove useful 
nowadays in domains as diverse as the laser control of quantum 
phenomena (see~\cite{hbref3,hbref10,hsiehscience,hbref8,hbref9,hbref4}) 
or the modeling of a equilibrium (or again social beliefs, 
product prices, etc) of a game with an infinite numbers of agents (see~\cite{mfg1,mfg2,mfg3}). Yet other
applications arise from modern formulations of the 
{M}onge-{K}antorovich mass transfer problem (see~\cite{MR1738163,MR1865668,cdcCarlierSalomon}). 

 In the following, we present some examples coming from these fields of application
and present the corresponding monotonic algorithm resulting from Theorem~\ref{thm:existence}.

\subsection{(I): Quantum control} \label{sec:exqc}
\subsubsection{Setting}
The evolution of a quantum system is 
described by the 
Schr\"odinger equation
\beqn \nonumber
& \ & \partial_t X  + i H(t) X = 0 \\
& \ & \nonumber X(0,z) = X_0(z),
\eeqn
where $i=\sqrt{-1}$,  $H(t)$ is the Hamiltonian of the system and $z \in \R^\gamma$ the set of
internal degrees of freedom. We assume that the Hamiltonian is a self-adjoint operator
over ${L^2({\R}^{{\gamma}}; \C)}$, i.e.
$H(t)^* = H(t)$\footnote{For any operator $M$,
we denote by $M^*$ its adjoint.}.
Note that this implies the following norm conservation
property
\beq \label{eq:consv} 
\|X(t,\cdot)\|_{L^2({\R}^{{\gamma}}; \C)}=
\|X_0\|_{L^2({\R}^{{\gamma}}; \C)},
 \ \forall t>0,
\eeq
so that the state (also called wave-) function $X(t,z)$, evolves
on the (complex) unit sphere $S \triangleq \left\{ {X  \in {L^2({\R}^{{\gamma}}; \C)}\,:
\,\left\| X  \right\|_{L^2({\R}^{{\gamma}}; \C)}  = 1} \right\}$.

The Hamiltonian is composed of two parts: a free evolution Hamiltonian $H_0$ and
a part that describes the coupling of the system with an
external laser source of intensity
$v(t) \in \R, \ t \ge 0$; a first order approximation
leads to adding a time-independent dipole moment operator $\mu(z)$ resulting
in the formula $H(t) = H_0 -v(t)\mu$ and the dynamics:
\beqn \nonumber 
& \ & \partial_t X + i\left( H_0 -v(t)\mu \right) X=0  \\
& \ & \nonumber X(0) = X_0.
\eeqn

The purpose of control may be formulated as
to drive the system from its initial state $X_0$ 
to a final state $X_{target}$ compatible with
predefined requirements. Here, the control is the laser intensity $v(t)$.
Because the control is multiplying the state, this formulation is called
``bilinear'' control. The dependence $v \mapsto X(T)$ is of course not linear.

The optimal control approach can be implemented by introducing a cost
functional. The following functionals are often
  considered:
\beqn \label{fonctquant0}
& \ & 
J(v) \triangleq \| X(T) - X_{target}\|^2_{L^2({\R}^{{\gamma}}; \C)} + \int_0^T \alpha(t) v^2 (t) dt,
\\ & \ & 
\widetilde{J}(v) \triangleq -\langle X(T),O X(T)\rangle_{L^2({\R}^{{\gamma}}; \C)} + \int_0^T \alpha(t) v^2 (t) dt,
\eeqn
where $O$ is a positive linear operator defined on $\Hspace$,
  characterizing an observable quantity and
$\alpha(t) > 0$ is a parameter that penalizes large (in the $L^2$ 
sense) controls. The goal is
to minimize these functionals with respect to $v$. According to~\eqref{eq:consv}
the cost functional $J$ is equal to
\beq \label{fonctquant}
J(v) \triangleq 2 - 2 Re \langle X(T), X_{target}\rangle_{L^2({\R}^{{\gamma}}; \C)} + \int_0^T \alpha(t) v^2 (t) dt,
\eeq
so that 
the functionals $J$ and $\widetilde{J}$ satisfy assumptions~\eqref{concavG} and \eqref{concavF}.

\subsubsection{Mathematical formulation}
We have 

$\bullet$ $A(t,v) = H_0 + v(t) \mu $ with (possibly) unbounded $v$-independent operator $H_0$ (but which 
generates a $C^0$ semi group) and bounded operator $\mu$. The dependence of $A$ on $v$ is smooth (linear) and therefore all
hypotheses on $A$ are satisfied. 

$\bullet$  $E=\R$, 
$\Hspace=L^2(\R^d;\C)$, $\Vspace=  dom(H_0^{1/2})$ (over $\C$), or their realifications 
$\Hspace=L^2 \times L^2$, $\Vspace=  dom(H_0^{1/2}) \times dom(H_0^{1/2})$ (over $\R$) as explained in~\cite{kunish2};

$\bullet$ $B(t,v)=0$.

$\bullet$ $F(t,v,X)=\alpha(t) v(t)^2$ with $\alpha(t)\in L^\infty(\R)$; here the second derivative $D_{vv}F$ is obviously bounded. Since it is
independent of $X$ it will be trivially concave.

$\bullet$ $G$ is either (see e.g.,~\cite{maday:2468,jcp}) $2 -2Re\langle X_{target},X(T)\rangle_{\Vspace}$ or 
$-\!\langle X(T),OX(T)\rangle_{\Vspace}$ where $O$ is a positive semi-definite operator; both are concave in $X$.

$\bullet$ Here 
\beq 
\Delta(v',v;t,X,Y) = - Re \langle Y , i\mu X\rangle_{\Vspace} + \alpha(t)(v'+v)
\eeq
and the equation in $v'$ is: 
$\Delta(v',v;t,X,Y) = - \theta  (v'-v)$ and has for $\theta$ large enough a unique solution
$v'={\cal V}_\theta(t,v,X,Y) \triangleq \frac{ (\theta-\alpha(t)) v +  Re \langle Y , i\mu X\rangle_{\Vspace}}{\theta + \alpha(t)}$.

$\bullet$ at the $k+1$-th iteration, Theorem~\ref{thm:existence} guarantees the existence of the solution
$X^{k+1}$ of the following nonlinear evolution equation:
\beq
i\partial_t{X^{k+1}}(t)= \left( H_0 + \frac{ (\theta-\alpha(t)) v^k +  Re \langle Y_{v^{k}} , i\mu X^{k+1}\rangle_{\Vspace}}{\theta + \alpha(t)} \mu \right) X^{k+1}(t) \label{eq:Xnl12}
\eeq
Then 
\beq
v^{k+1}= \frac{ (\theta-\alpha(t)) v^k + Re \langle Y_{v^{k}} , i\mu X^{k+1}\rangle_{\Vspace}}{\theta + \alpha(t)}, \ 
X_{v^{k+1}}=X^{k+1}.\label{eq:itqc}
\eeq
\subsubsection{Numerical test} \label{sec:exqcnum}
In order to test the performance of the algorithm we have chosen a
case already treated in the literature (see~\cite{zr98}). 
The
system under consideration is the $O-H$ bond that vibrates in a Morse
type potential $V(z) = D_0(\exp(-\beta(z- z')) -1 )^2 - D_0$ and
$H_0=-m\frac{\partial^2}{\partial z^2}+V(z)$. The
dipole moment operator of this system is modeled by $\mu(z)=\mu_0.z e^{-\frac z{ z^\star}}$.
 The objective is to localize the wavefunction at time $T=131000$ at a given
location $z_0$ ;
this is expressed through the requirement that the functional
$\widetilde{J}$ is maximized, with the observable $O$ defined 
by $O(z) = \frac{\gamma_0}{\sqrt{\pi}}e^{-\gamma_0^2(z-z_0)^2}$. 
The numerical values we use are given below:
$$
\begin{array}{|c|c|c|c|c|c|c|c|}
\hline
D_0    & \beta & z'     & z^\star & z_0    & \gamma_0 & \mu_0 & m\\\hline
0.1994 & 1.189 &  1.821 & 0.6     & 2.5    &       25 & 3.088 & 2.8694.10^{-4}\\ \hline
\end{array}
$$

\begin{figure}[h]
\begin{center}
\psfrag{a}[c][t]{iterations}
\psfrag{b}[c][t]{$J(v^k)$}
  \includegraphics[width=.46\textwidth]{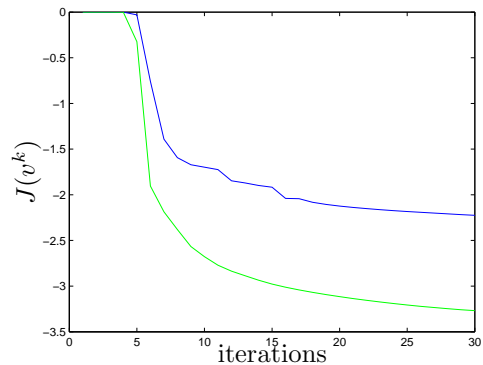}
\caption{Numerical resolution of the example of
  Section~\ref{sec:exqc}. The cost functional $J(v)$ is optimized
  using the
monotonic algorithm~\eqref{eq:itqc} (green line) and the optimal step gradient algorithm
(blue line). }\label{ex1}
\end{center}
\end{figure}

We
consider a constant
penalization parameter $\alpha=1$ and optimization parameter $\theta=10^{-2}$. 
To compare this procedure with a standard algorithm, we have also
minimized $J(v)$
with an optimal step gradient method. The line
search is achieved through a golden section search cf.~\cite{nr}. 
Results are presented in Fig.~\ref{ex1}. 
This test shows that the gradient method fails in efficiently solving 
the problem, whereas the monotonic procedure ensures that the cost functional
values rapidly decrease. Note that the non-convexity of the problem
renders difficult the convergence of the gradient method. On the other hand, the
monotonic scheme fully exploits the concavity of the cost functional
with respect to the state 
variable. In our implementation the time of
computation is about two times larger for the gradient method as the
line search requires about 3 evaluations of the cost functional per
iteration.

\subsection{(II) : Mean field games} \label{sec:exmfg}

\subsubsection{Setting}

Although the Nash equilibrium in game theory has been initially formulated for a finite number of players, modern
results (see~\cite{mfg1,mfg2,mfg3}) indicate that it is possible to extend 
it to an infinite number of players and obtain the equations that describe this equilibrium; applications 
have already been proposed in economic theory and other are expected in the behavior of multi-agents ensembles
and decision theory. 

The equations describe evolution of the density $X(t,z)$ of players at time $t$ and position $z \in Q=[0,1]$ in terms
of a control $v(t,{z})$ and a fixed parameter $\nu >0$: 
\beqn \nonumber 
& \ & \partial_t X  - \nu \Delta X + div(v(t,{z})X) = 0,\\
& \ & \nonumber X(0) = X_0.
\eeqn

The control $v$ is chosen to minimize the cost criterion~\eqref{eqn:F}. For reasons related to economic modeling
interesting examples include situations where $F,G$ are concave in $X$, 
e.g., as in~\cite{aimemfg}
\beq \label{eq:L}
G=0, \  F(t,z,X)=  \int_Q  p(t)(1-\beta z)X(t,z) + \frac{c_0 \cdot z
  \cdot X(t,z)}{c_1+c_2 X(t,z)} + \frac{v^2(t)}{2}X(t,z)dz, 
\eeq
with positive constants $\beta,c_0,c_1,c_2$ and 
$p(t)$ a positive function. Another example is given in~\cite{cdcCarlierSalomon}:
\beq \label{eq:LCDC}
G(X(T))=\int_Q V(z) X(T,z) dz, \  F(t,z,X)=\int_Q X(t,z) v^2(t,z) dz,
\eeq
where $V$ encodes a potential. The interpretation of this terminal
cost is that the crowd aims at reaching zones of low potential $V$ at the  
terminal time $T$ while minimizing the cost of changing state.

The numerical relevance of the monotonic algorithms to this setting has been
established in several works, see~\cite{cdcCarlierSalomon,aimemfg}.

\subsubsection{Mathematical formulation}

We have 

$\bullet$ $E=W^{1,\infty}(0,1)$,  $\Hspace=L^2(0,1)$, $\Vspace=  H^1(0,1) $ see~\cite{aimemfg} and~\cite{dautraylions5} (Chap XVIII $\S$4.4)

$\bullet$ $A(t,v) = -\nu \Delta \cdot + div(v \cdot) $. The dependence of $A$ on $v$ is smooth (linear) and therefore all
hypotheses on $A$ are satisfied ($D_{vv}A=0, ...$). 

$\bullet$ $B(t,v)=0$.

$\bullet$  with definitions in~\eqref{eq:L} 
$F(t,v,X)=\int_Q p(t)(1-\beta z)X(t,z) + \frac{c_0 \cdot z \cdot X(t,z)}{c_1+c_2X(t,z)}  + \frac{v(t,z)^2}{2} X(t,z)dz$; $F$ is concave 
in $X$; 
the second differential $D_{vv}F$ has all required properties.

$\bullet$ $G=0$ (algorithm will apply in general when $G$ is concave with respect to $X$).

$\bullet$ Here
\beq
\Delta(v',v;t,X,Y) = \nabla Y + \frac{v'+v}{2}
\eeq
and the equation in $v'$ is: 
 $\Delta(v',v;t,X,Y) = - \theta  (v'-v)$
and has for all $\theta>0$ a unique solution
$v'={\cal V}_\theta(t,v,X,Y) \triangleq \frac{ (\theta-1/2) v -\nabla Y}{\theta + 1/2}$.

$\bullet$ at the $k+1$-th iteration, Theorem~\ref{thm:existence} guarantees the existence of the solution
$X^{k+1}$ of the following nonlinear evolution equation:
\beq
\partial_t{X^{k+1}}(t) -\nu \Delta X^{k+1} + div(\frac{ (\theta-1/2)
  v^k -\nabla Y_{v^k}}{\theta + 1/2} X^{k+1})  =0.
\eeq
Then 
\beq\label{eq:itmfg}
v^{k+1}=\frac{ (\theta-1/2) v^k -\nabla Y_{v^k}}{\theta + 1/2}, \ X_{v^{k+1}}=X^{k+1}.
\eeq

\subsubsection{Numerical test}
\begin{figure}[h]
\begin{center}
\psfrag{a}[c][t]{iterations}
\psfrag{b}[c][t]{$J(v^k)$}
  \includegraphics[width=.46\textwidth]{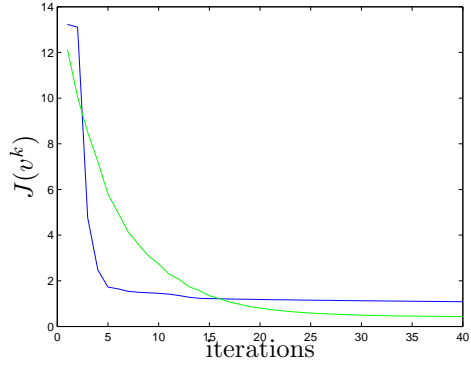}
\caption{Numerical resolution of the example of
  Section~\ref{sec:exmfg}. The cost functional $J(v)$ is optimized
  using the
monotonic algorithm~\eqref{eq:itmfg} (green line) 
and the optimal step gradient algorithm (blue line). } \label{ex2}
\end{center}
\end{figure}

The algorithm is tested on the time interval $[0,1]$ with $p(t)=1$ and
the numerical values $\beta=0.8$, $c_0= c_2=1$ , $c_1=0.1$. The same gradient
method as in Section~\ref{sec:exqcnum} is also tested. Results are presented in Fig.~\ref{ex2}. 
In this example, the gradient method gives better results in the first
iterations. However, the monotonic algorithm converges asymptotically
faster.

\subsection{Additional application} \label{sec:exnonlin}

As a third example we consider a 
nonlinear vectorial case from~\cite{friedrich,tehini} which differs from that of 
Section~\ref{sec:exqc} in that $v(t)=
\begin{pmatrix}
v_1\\ v_2 
\end{pmatrix}
 \in E=\R^2$ and $A(t,v)= i [ H_0 + (v_1(t)^2+ v_2(t)^2) \mu_1 +  v_1(t)^2 v_2(t) \mu_2 ]$.
Here, denoting
$\xi_1 = 
-Re \langle Y , i \mu_1 X\rangle_{\Vspace} + \alpha(t) 
$,
$\xi_2= -Re \langle Y , i \mu_2 X\rangle_{\Vspace}$ 
we obtain 
\beq\label{eq:itnonlin}
\Delta(v',v;t,X,Y) = \xi_1 \begin{pmatrix}
v_1 +v'_1 \\ v_2 +v'_2  
\end{pmatrix}
+\xi_2 \begin{pmatrix}
(v_1 +v'_1)v'_2 \\ (v_1)^2 
\end{pmatrix}
\eeq
and the equation in $v'$ is: 
$\Delta(v',v;t,X,Y) = - \theta  (v'-v)$
and has for $\theta$ large enough  a unique solution
$v'={\cal V}_\theta(t,v,X,Y)= 
\begin{pmatrix}
\frac{ (\theta - \xi_1) v_2 - \xi_2 v_1^2}{\theta + \xi_1} \\
- \frac{ \theta - \xi_1 + \xi_2  \frac{ (\theta - \xi_1) v_2 - \xi_2 v_1^2}{\theta + \xi_1} }{  \theta + \xi_1 + \xi_2  \frac{ (\theta - \xi_1) v_2 - \xi_2 v_1^2}{\theta + \xi_1}}v_1
\end{pmatrix}$.
We leave as an exercise to the reader the writing of the equation for $X^{k+1}$ and the formula for $v^{k+1}$.\\
This model corresponds to the problem of controlling the orientation $\gamma$ of
a molecule, considered as rigid rotator.

\subsubsection{Numerical test}

\begin{figure}[h]
\begin{center}
\psfrag{a}[c][t]{iterations}
\psfrag{b}[c][t]{$J(v^k)$}
  \includegraphics[width=.46\textwidth]{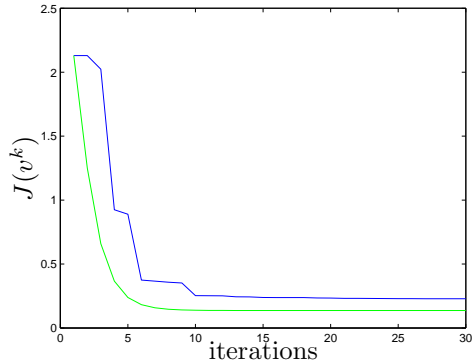}
\caption{Numerical resolution of the example of
  Section~\ref{sec:exnonlin}. The cost functional $J(v)$ is optimized
  using the
monotonic algorithm~\eqref{eq:itnonlin} (green line) and a gradient method
(blue line).}\label{ex3}\end{center}
\end{figure}

To test our approach we have used the parameters of the molecule
$CO$ (see~\cite{friedrich,tehini}), namely $H_0= BJ^2$, where $B$ is the 
rotational constant and $J$ is the angular momentum. We consider the
basis given by the spherical harmonics~; the corresponding matrix is
diagonal with diagonal coefficients given by $(H_0)_{k,k}=k(k+1)$.
The controlled is performed over an interval of length
$T=20T_{per}=20\frac\pi B$. We consider constant penalization factor
$\alpha=10^{-1}$ and optimization parameter $\theta=10^{3}$.\\
The other parameters correspond to the polarizability and the
hyperpolarizability components of the molecule. We have
$\mu_1= -\frac12\lambda,$ and $\mu_2= -\frac34\beta,$ with
$\lambda=\frac12(\lambda_\para\cos^2\gamma+\lambda_\perp\sin^2\gamma),$ 
 $\beta=\frac16((\beta_\para-3\beta_\perp)\cos^3\gamma+3\beta_\perp\cos\gamma).$
The matrix $\cos \gamma$ is tridiagonal, with:
$$ (\cos \gamma)_{k,k}=0, \  (\cos \gamma)_{k,k+1}=(\cos \gamma)_{k+1,k}= \frac{k+1}{\sqrt{(2k+1)(2k+3)}}.$$
We use the numerical values given in~\cite{friedrich,tehini}:  
$$
\begin{array}{|c|c|c|c|c|}
\hline
B & \lambda_\perp &  \lambda_\para & \beta_\para & \beta_\perp \\
\hline
1.93 & 11.73 & 15.65 & 28.35 & 6.64\\ \hline
\end{array}
$$
 A gradient
method similar to the one that is used in Section~\ref{sec:exqcnum} is
also performed.
The results are presented in Fig.~\ref{ex3}. The monotonic algorithm
shows a fast convergence whereas the gradient method does not
optimizes efficiently the cost functional values.

\section{Conclusion}

Motivated by a set of control algorithms that were initially
introduced in the specific context of quantum control
we have presented an abstract formulation that includes them all.
It is seen that the algorithm involves at each step the resolution of a highly nonlinear evolution equation.
We identified the theoretical assumptions that ensure that the evolution equation
is well posed and has a solution. The proof being constructive it 
serves as basis for numerical approximations of the solution.
We also proved several properties concerning the algorithms and more specifically
concerning its convergence.
Examples are provided to indicate how the proposed procedure 
solves cases from the literature and also new
situations that were not previously considered. Numerical simulations indicate that the procedures 
have indeed the expected behavior.

\section*{Acknowledgements}
This work is partially supported by the French ANR programs
 OTARIE (ANR-07-BLAN-0235 OTARIE) C-QUID (grant BLAN-3-139579) and by a CNRS-NFS PICS grant. G.T. acknowledges partial support by INRIA Rocquencourt (MicMac and OMQP).

\bibliographystyle{plain}
\bibliography{references}
\vspace{12pt}

\end{document}